\documentclass{IEEEtran}
\usepackage{cite}
\usepackage{amsmath,amssymb,amsfonts}
\usepackage{algorithmic}
\usepackage{graphicx}
\usepackage{textcomp}
\def\BibTeX{{\rm B\kern-.05em{\sc i\kern-.025em b}\kern-.08em
    T\kern-.1667em\lower.7ex\hbox{E}\kern-.125emX}}
\usepackage{enumerate}
\usepackage{enumitem}
\usepackage{mathrsfs}
\usepackage{cases}
\usepackage{color}
\usepackage{tikz}
\usepackage{pgfplots}
\pgfplotsset{width=14cm,height=6cm,compat=1.9}
\usepackage{makecell}
\usepackage{multirow}
\usepackage{subfig}
\definecolor{darkblue}{rgb}{0.0, 0.0, 0.55}
\definecolor{bordeaux}{rgb}{0.34, 0.01, 0.1}
\usepackage[colorlinks,linkcolor=bordeaux,citecolor=darkblue,urlcolor=black,hypertexnames=true]{hyperref}
\usetikzlibrary{arrows,positioning}
\usepackage{mathtools}
\definecolor{pink1}{RGB}{255,181,197}
\definecolor{lightblue}{RGB}{173,216,230}
\definecolor{lightgreen}{RGB}{144,238,144}

\newtheorem{theorem}{Theorem}[section]
\newtheorem{lemma}[theorem]{Lemma}
\newtheorem{definition}[theorem]{Definition}
\newtheorem{example}[theorem]{Example}

\newtheorem{proposition}[theorem]{Proposition}

\newtheorem{remark}[theorem]{Remark}

\definecolor{dkgreen}{rgb}{0,0.4,0}

\def\Z{{\mathbb{Z}}}

\def\R{{\mathbb{R}}}

\def\N{{\mathbb{N}}}

\def\x{{\mathbf{x}}}

\def\ba{{\mathbf{a}}}
\def\bb{{\mathbf{b}}}

\def\br{{\mathbf{r}}}

\def\a{{\boldsymbol{\alpha}}}

\def\b{{\boldsymbol{\beta}}}
\def\g{{\boldsymbol{\gamma}}}

\def\A{{\mathscr{A}}}
\def\B{{\mathscr{B}}}

\def\S{{\mathbf{S}}}

\def\supp{\hbox{\rm{supp}}}

\def\BibTeX{{\rm B\kern-.05em{\sc i\kern-.025em b}\kern-.08em
    T\kern-.1667em\lower.7ex\hbox{E}\kern-.125emX}}

\usepackage{cancel}

\begin{document}

\title{Exploiting Term Sparsity in Moment-SOS hierarchy for Dynamical Systems}

\author{Jie Wang, Corbinian Schlosser, Milan Korda, and Victor Magron
\thanks{Manuscript submitted 18 November, 2021. This work was supported in part by the NSFC under grant 12171324. C. Schlosser, M. Korda, and V. Magron have been supported by European Union's Horizon 2020 research and innovation programme under the Marie Sk\l{}odowska-Curie Actions, grant agreement 813211 (POEMA).}
\thanks{J. Wang is with Academy of Mathematics and Systems Science, CAS, China (e-mail: wangjie212@amss.ac.cn).}
\thanks{C. Schlosser, M. Korda, V. Magron are with Laboratory for Analysis and Architecture of Systems, CNRS, France (e-mail: corbinian.schlosser@gmail.com, korda@laas.fr, vmagron@laas.fr).}}

\maketitle

\begin{abstract}
In this paper, we develop a dynamical system counterpart to the term sparsity sum-of-squares (TSSOS) algorithm proposed for static polynomial optimization. This allows for computational savings and improved scalability while preserving convergence guarantees when sum-of-squares methods are applied to problems from dynamical systems, including the problems of approximating region of attraction, the maximum positively invariant set, and the global attractor. At its core, the method exploits the algebraic structure of the data, thereby complementing existing methods that exploit causality relations among the states of the dynamical system. The procedure encompasses sign symmetries of the dynamical system as was already revealed for polynomial optimization. Numerical examples demonstrate the efficiency of the approach in the presence of this type of sparsity.
\end{abstract}


\section{Introduction}\label{intro}

The idea of translating problems from dynamical systems to infinite dimensional linear programming problems dates back to at least the work of Rubio~\cite{Rubio} and Lewis and Vinter~\cite{lewis1980relaxation} concerned with optimal control problems. More recently this idea was extended to other problems: approximations of maximum positively invariant (MPI) set \cite{korda2014convex}, region of attraction (ROA) \cite{kordaAttraction}, reachable set \cite{magron2019semidefinite}, global attractors (GA) \cite{MilanCorbiAttractor} and invariant measures \cite{MilanInvariant} among others. These problems then can be solved in the spirit of \cite{LasserreOCP} using a convergent sequence of finite dimensional convex optimization problems. This procedure results in a hierarchy of moment-sum-of-squares (moment-SOS) relaxations leading to a sequence of semidefinite programs (SDPs), as was done for optimal control problems in \cite{LasserreOCP}.

However, the size of these SDPs scales rapidly with the relaxation order and the state-space dimension. As a consequence, despite being convex, these SDP relaxations may be challenging to solve even for problems of modest state-space dimension. In the context of polynomial optimization, the similar problem of scalability has been extensively studied in recent years via exploiting structure of the system, e.g., \cite{riener2013exploiting} by exploiting symmetries, \cite{waki} by exploiting correlative sparsity, \cite{chordaltssos,tssos} by exploiting term sparsity.
In this paper we present the use of the recent term sparsity approach which has been already proven useful for a wide range of polynomial optimization problems, involving 
noncommutating variables \cite{nctssos}, and fast approximation of joint spectral radius of sparse matrices \cite{wang2020sparsejsr}.
For these problems, one is able to formulate computationally cheaper hierarchies by exploiting term sparsity with strong convergence properties. 

Whereas the approaches of \cite{schlosser2020sparse,tacchi2020approximating} are concerned with the sparsity in the couplings between variables themselves, the approach proposed here exploits sparsity in the \emph{algebraic description} of the dynamics, in particular among the monomial terms appearing in the components of the polynomial vector field. The method proceeds by searching for non-negativity certificates comprised of polynomials with only specific sets of terms, which in turn are enlarged in an iterative scheme. 
From an operator-theoretic perspective, the proposed term sparsity approach exploits term sparsity of the data and dynamics by algebraic (or graph theoretic) properties of the Liouville-operator associated with the dynamics. Interestingly, this approach intrinsically comprises the sign symmetry reduction. In general, the term sparsity approach allows a trade-off between the computational cost and the accuracy of approximation, whereas exploiting sign symmetries does not sacrifice any accuracy. So the term sparsity approach would provide an alternative reduction when no non-trivial sign symmetry is available or exploiting sign symmetries is still too expensive.

To summarize, the main contributions of this paper are: (1) We develop a term sparsity approach for SOS methods in dynamical systems, allowing for computational reduction beyond existing sparsity exploitation methods.
(2) The method provides a sequence of SDPs of increasing complexity and we show that this sequence converges in finitely many steps to the sign-symmetry reduction of the original problem.
(3) We demonstrate the approach on a number of examples (including the extended Lorenz system, randomly generated instances, and a 16-state fluid mechanics example) and observe a promising trade-off between speed-up and the solution accuracy.


The rest of this paper is organized as follows. In Section~\ref{pre}, we introduce the notation and give some preliminaries. In Section~\ref{ts}, we show how to exploit term sparsity in the moment-SOS hierarchy by taking the computation of MPI sets as an example and revealing its relation with the sign symmetry reduction. Section~\ref{example} illustrates the approach by numerical examples. Conclusions are given in Section~\ref{con}.


\section{Notation and preliminaries}\label{pre}

Let $\x=(x_1,\ldots,x_n)$ be a tuple of variables and $\R[\x]=\R[x_1,\ldots,x_n]$ be the ring of real $n$-variate polynomials. For $d\in\N$, the subset of polynomials in $\R[\x]$ of degree no more than $2d$ is denoted by $\R_{2d}[\x]$. A polynomial $f\in\R[\x]$ can be written as $f(\x)=\sum_{\a\in\A}f_{\a}\x^{\a}$ with $\A\subseteq\N^n$, $f_{\a}\in\R$ and $\x^{\a}=x_1^{\alpha_1}\cdots x_n^{\alpha_n}$. The support of $f$ is then defined by $\supp(f)\coloneqq\{\a\in\A\mid f_{\a}\ne0\}$. For $\A\subseteq\N^n$, let $\R[\A]$ be the set of polynomials whose supports are contained in $\A$. The notation $Q\succeq0$ for a matrix $Q$ indicates that $Q$ is positive semidefinite (PSD). For a positive integer $r$, the set of $r\times r$ symmetric matrices is denoted by $\mathbf{S}^r$, and the set of $r\times r$ PSD matrices is denoted by $\mathbf{S}_+^r$. For $d\in\N$, let $\N^n_d\coloneqq\{\a=(\alpha_i)_{i=1}^n\in\N^n\mid\sum_{i=1}^n\alpha_i\le d\}$. For $\a\in\N^n,\A,\B\subseteq\N^n$, let $\a+\B\coloneqq\{\a+\b\mid\b\in\B\}$ and $\A+\B\coloneqq\{\a+\b\mid\a\in\A,\b\in\B\}$. We use $|\cdot|$ to denote the cardinality of a set.
For two vectors $\ba=(a_i)_{i=1}^n$ and $\bb=(b_i)_{i=1}^n$, let $\ba\cdot\bb\coloneqq\sum_{i=1}^na_ib_i$ and $\ba\circ\bb\coloneqq(a_1b_1,\ldots,a_nb_n)$.

Given a polynomial $f(\x)\in\R[\x]$, if there exist polynomials $f_1(\x),\ldots,f_t(\x)$ such that $f(\x)=\sum_{i=1}^tf_i(\x)^2$,
then we call $f(\x)$ a {\em sum of squares (SOS)} polynomial. The set of SOS polynomials is denoted by $\Sigma[\x]$. Assume that $f\in\Sigma_{2d}[\x]\coloneqq\Sigma[\x]\cap\R_{2d}[\x]$ and $\x^{\N^n_{d}}$ is the $\binom{n+d}{d}$-dimensional column vector consisting of elements $\x^{\a},\a\in\N^n_{d}$ (fix any ordering on $\N^n$). Then $f$ is an SOS polynomial if and only if there exists a PSD matrix $Q$ (called a Gram matrix) such that $f=(\x^{\N^n_{d}})^\intercal Q\x^{\N^n_{d}}$.
For convenience, we abuse notation in the sequel and denote by $\N^n_{d}$ instead of $\x^{\N^n_{d}}$ the standard monomial basis and use the exponent $\a$ to represent a monomial $\x^{\a}$.

An {\em (undirected) graph} $G(V,E)$, or simply $G$, consists of a set of nodes $V$ and a set of edges $E\subseteq\{\{u,v\}\mid u\ne v,(u,v)\in V\times V\}$. For a graph $G$, we use $V(G)$ and $E(G)$ to indicate the node set of $G$ and the edge set of $G$, respectively. For two graphs $G,H$, we say that $G$ is a {\em subgraph} of $H$, denoted by $G\subseteq H$, if both $V(G)\subseteq V(H)$ and $E(G)\subseteq E(H)$ hold.
A graph is called a {\em chordal graph} if all its cycles of length at least four have a chord\footnote{A chord is an edge that joins two nonconsecutive nodes in a cycle.}.
The notion of chordal graphs plays an important role in sparse matrix theory. Any non-chordal graph $G(V,E)$ can be always extended to a chordal graph $G'(V,E')$ by adding appropriate edges to $E$, which is called a {\em chordal extension} of $G(V,E)$. 
The chordal extension of $G$ is usually not unique and the symbol $G'$ is used to represent any specific chordal extension of $G$ throughout the paper. For a graph $G$, there is a particular chordal extension which makes every connected component of $G$ to be a complete subgraph, which is called the {\em maximal} chordal extension. Typically, we consider only chordal extensions that are subgraphs of the maximal chordal extension. For graphs $G\subseteq H$, we assume that $G'\subseteq H'$ holds throughout the paper. 
Given a graph $G(V,E)$, a symmetric matrix $Q$ with rows and columns indexed by $V$ is said to have sparsity graph $G$ if $Q_{uv}=Q_{vu}=0$ whenever $u\ne v$ and $\{u,v\}\notin E$. Let $\mathbf{S}_G$ be the set of symmetric matrices with sparsity graph $G$. 
The PSD matrices with sparsity graph $G$ form a convex cone $\mathbf{S}_+^{|V|}\cap\mathbf{S}_G=\{Q\in\mathbf{S}_G\mid Q\succeq0\}$.
When the sparsity graph $G(V,E)$ is a chordal graph, the cone $\mathbf{S}_+^{|V|}\cap\mathbf{S}_G$ can be
decomposed as a sum of simple convex cones by virtue of the following theorem. Recall that a {\em clique} of a graph is a subset of nodes that induces a complete subgraph. A {\em maximal clique} is a clique that is not contained in any other clique. 
\begin{theorem}[\cite{agler}, Theorem 2.3]\label{sec2-thm}
Let $G(V,E)$ be a chordal graph and assume that $C_1,\ldots,C_t$ are all maximal cliques of $G(V,E)$. Then a matrix $Q\in\mathbf{S}_+^{|V|}\cap\mathbf{S}_G$ if and only if $Q$ can be written as $Q=\sum_{i=1}^tQ_{i}$, where $Q_i\in\mathbf{S}_+^{|V|}$ has nonzero entries only with row and column indices coming from $C_i$.
\end{theorem}
The chordal decomposition stated in Theorem \ref{sec2-thm} has enabled significant progress in large-scale semidefinite programming; see for instance \cite{fukuda2001exploiting,va}.
Given a graph $G$ with $V=\N^n_{d}$, let $\Sigma[G]$ be the set of SOS polynomials that admit a PSD Gram matrix with sparsity graph $G$, i.e., $\Sigma[G]\coloneqq\{(\x^{\N^n_{d}})^\intercal Q\x^{\N^n_{d}}\mid Q\in\mathbf{S}_+^{|V|}\cap\mathbf{S}_G\}$. Note that in general $\Sigma[G]$ is a strict subset of $\Sigma_{2d}[\x]$, and therefore the sparse SOS strengthening $f\in\Sigma[G]$ of the inequality $f\ge0$ is generally more conservative than the corresponding dense strengthening $f\in\Sigma_{2d}[\x]$.

\section{Exploiting term sparsity}\label{ts}
In this section, we propose an iterative procedure to exploit term sparsity for the moment-SOS hierarchy of certain computational problems related to dynamical systems. The intuition behind this procedure is the following: starting with a minimal initial support set, we expand the support set that is taken into account in the SOS relaxation by iteratively performing \emph{support extension} and \emph{chordal extension} to the related sparsity graphs inspired by Theorem \ref{sec2-thm}. In doing so, we obtain ascending chains of support sets which further lead to a hierarchy of sparse SDP relaxations.
For the ease of understanding, we illustrate the approach by considering the computation of MPI sets. But there is no difficulty to extend the approach to other situations, e.g., the computations of ROA \cite{kordaAttraction} and GA \cite{MilanCorbiAttractor}, bounding extreme events \cite{fantuzzi2020bounding}.
Suppose that the dynamical system we are considering is given by
\begin{equation}\label{dyn}
\dot{\x}=\mathbf{f}(\x),
\end{equation}
with $\mathbf{f}\coloneqq\{f_1,\ldots,f_n\}\subseteq\R[\x]$, and the constraint set is
\begin{equation}\label{cons}
X\coloneqq\{\x\in\R^n\mid p_j(\x)\ge0\textrm{ for }j=1,\ldots,m\}
\end{equation}
with $p_1,\ldots,p_m\in\R[\x]$. For the sake of convenience, we set $p_0\coloneqq1$. Let $d_f\coloneqq\max\,\{\deg(f_i):i=1,\ldots,n\}$ and $d_j\coloneqq\deg(p_j)$ for $j=0,1,\ldots,m$, $d_p\coloneqq\max\,\{d_j:j=1,\ldots,m\}$. For $\x_0\in\R^n$ and $t\ge0$, let $\varphi_t(\x_0)$ denote the solution
of \eqref{dyn} with initial condition $\varphi_0(\x_0)=\x_0$.
\begin{definition}
For a dynamical system \eqref{dyn} with the constraint set $X$, the {\em maximum positively invariant (MPI) set} is the set of initial conditions $\x_0\in X$ such that the solutions $\varphi_t(\x_0)$ stay in $X$ for all $t\ge0$.
\end{definition}

As proposed in \cite{korda2014convex}, given a positive integer $d$, the $d$-th order SOS relaxation for approximating the MPI set is defined by
\begin{equation}\label{sos}
\theta_{d}\coloneqq\begin{cases}    \inf\limits_{a_j,b_j,c_j,v,w}\quad&\int_{X}w(\x)\,\mathrm{d}\x\\
    \quad\,\,\,\text{s.t.}&v\in\R[\x]_{2d+1-d_f},w\in\R[\x]_{2d},\\
    &\beta v-\nabla v\cdot \mathbf{f}=a_0+\sum_{j=1}^ma_jp_j,\\
    &w=b_0+\sum_{j=1}^mb_jp_j,\\
    &w-v-1=c_0+\sum_{j=1}^mc_jp_j,\\
    &a_j,b_j,c_j\in\Sigma_{2d-d_j}[\x],j=0,1,\ldots,m,
\end{cases}
\end{equation}
where $\beta>0$ is a preassigned discount factor (say, $\beta=1$), $\nabla$ is the gradient with respect to $\x$. The dynamics enter through the discounted Liouville operator $v \mapsto \beta v- \nabla v \cdot \mathbf{f}$. By \cite{korda2014convex}, the sequence of optima of \eqref{sos} converges monotonically from above to the volume of the MPI set provided the polynomials $(p_j)_{j=1}^m$ satisfy a technical compactness condition, which is satisfied for example if $p_m=R-\|\x\|^2$ for some $R>0$.
Furthermore, the set
$S_d\coloneqq w^{-1}([1,+\infty])=\{\x\in X: w(\x)\ge1\}$
provides an outer approximation for the MPI set.

For a graph $G(V,E)$ with $V\subseteq\N^n$, define $\supp(G)\coloneqq\{\b+\g\mid\b=\g\in V\text{ or }\{\b,\g\}\in E\}$. Now we give the iterative procedure to exploit term sparsity. Fix a relaxation order $d$. Let $\A=\bigcup_{j=1}^m\supp(p_j)$ and $v$ be a polynomial with generic coefficients supported on $\A$. Let 
$\A_{d}^1\coloneqq\A\cup\supp(\nabla v\cdot \mathbf{f})\cup2\N^n_d$
with $2\N^n_d\coloneqq\{2\a\mid\a\in\N^n_d\}$. Intuitively, the set $\A_{d}^1$ is the minimal support that has to be involved in the SOS relaxation \eqref{sos}\footnote{Here the subset $2\N^n_d$ is included in the definition of $\A_{d}^1$ to guarantee convergence; see \cite{tssos}.}. Now we expand the support set by iteratively performing support extension and chordal extension on some related graphs. In the following we perform a construction which is motivated by the desire to impose sparsity patterns to Gram matrices of $\beta v-\nabla v\cdot \mathbf{f},w$ and $w-v-1$ which are required to be SOS polynomials over $X$ with a prescribed support. We begin with the construction for $\beta v-\nabla v\cdot \mathbf{f}$. For every integer $s\ge1$, we iteratively define the graph $G_{d,j}^s$ which will be imposed as the sparsity graph for a Gram matrix of $a_j$ with $V(G_{d,j}^s)\coloneqq\N^n_{d-d_j}$ and
\begin{equation}\label{sec3:eq1}
\begin{split}
E(G_{d,j}^s)\coloneqq\{\{\b,\g\}\mid\,&\b+\g+\supp(p_j)\\
&\cap(\A_{d}^s\cup\supp(\nabla v_{d}^s\cdot \mathbf{f}))\ne\emptyset\},
\end{split}
\end{equation}
where $v_d^s$ is a polynomial with generic coefficients supported on $\A_d^s\cap\N^n_{2d+1-d_f}$ for $j=0,1,\ldots,m$, and further let $\A_d^{s+1}=\supp((G_{d,0}^{s})')$.
Doing so, we get a finite ascending chain of support sets
$\A_d^{1}\subseteq\cdots\subseteq\A_d^{\tilde{s}}=\A_d^{\tilde{s}+1}=\cdots$
and a finite ascending chain of graphs
$(G_{d,j}^{1})'\subseteq\cdots\subseteq(G_{d,j}^{\tilde{s}})'=(G_{d,j}^{\tilde{s}+1})'=\cdots$
for each $j=0,1,\ldots,m$.

\begin{remark}
There is term sparsity to exploit for the SOS relaxation \eqref{sos} if the graph $G_{d,0}^{1}$ is not complete.
\end{remark}

\begin{example}\label{sec3:ex1}
Let us consider the classical Lorenz system: $\{
    \dot{x_1}=10(x_2-x_1),
    \dot{x_2}=x_1(28-x_3)-x_2,
    \dot{x_3}=x_1x_2-\frac{8}{3}x_3\}$
with the constraint set $X=\{\x\in\R^3\mid 1-x_1^2\ge0,1-x_2^2\ge0,1-x_3^2\ge0\}.$ Take the relaxation order $d=2$ and then $\A_{2}^1=\{1,x_1^2,x_2^2,x_3^2,x_1x_2,x_1x_2x_3,x_1^2x_2^2,x_1^2x_3^2,x_2^2x_3^2,x_1^4,x_2^4,x_3^4\}$. The sparsity graph $G_{2,0}^{1}$ is displayed in Fig.~\ref{supp1}. If we use maximal chordal extensions, then $(G_{2,0}^{1})'$ consists of two cliques respectively of size $6,4$, and $\A_{2}^2=\{1,x_1^2,x_2^2,x_3^2,x_1x_2,x_1x_2x_3,x_1^2x_2^2,x_1^2x_3^2,x_2^2x_3^2,x_1^4,x_2^4,x_3^4,x_3,$ $x_3^3,x_1x_2x_3^2,x_1^3x_2,x_1x_2^3\}$. In fact, we have $\A_{2}^2=\A_{2}^3=\cdots$ and $(G_{2,j}^{2})'=(G_{2,j}^{3})'=\cdots$ for $j=0,1,2,3$.

\begin{figure}[htbp]
\centering
{\tiny
\begin{tikzpicture}[scale=0.8, every node/.style={circle, draw=blue!50, thick, minimum size=7mm}]
\node (n1) at (-2,0) {$x_3$};
\node (n2) at (0,0) {$x_1$};
\node (n3) at (2,0) {$x_2$};
\node (n4) at (-2,-2) {$x_1x_2$};
\node (n5) at (0,-2) {$x_2x_3$};
\node (n6) at (2,-2) {$x_1x_3$};
\node (n7) at (-6,0) {$1$};
\node (n8) at (-4,0) {$x_1^2$};
\node (n9) at (-6,-2) {$x_3^2$};
\node (n10) at (-4,-2) {$x_2^2$};
\draw (n1)--(n4);
\draw (n7)--(n8);
\draw (n7)--(n9);
\draw (n7)--(n10);
\draw (n8)--(n9);
\draw (n8)--(n10);
\draw (n9)--(n10);
\draw (n2)--(n5);
\draw (n3)--(n6);
\draw (n1)--(n8);
\draw (n1)--(n10);
\draw (n2)--(n3);
\draw (n2)--(n6);
\draw (n3)--(n5);
\draw (n4)--(n7);
\end{tikzpicture}}
\caption{\small The sparsity graph $G_{2,0}^{1}$ for Example \ref{sec3:ex1}.}\label{supp1}
\end{figure}
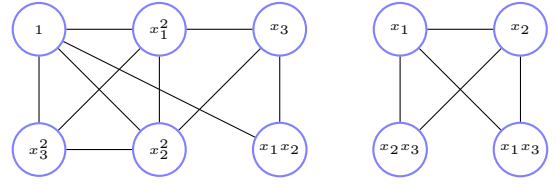

\end{example}

Now we perform similar constructions for $w$ and $w-v-1$.
For a given $s\ge1$, with $\B_d^{s,1}\coloneqq\A_d^s$ for every $l\ge1$ we iteratively define the graph $H_{d,j}^{s,l}$ which will be imposed as the sparsity graph for a Gram matrix of $b_j$ or $c_j$ with $V(H_{d,j}^{s,l})\coloneqq\N^n_{d-d_j}$ and
\begin{equation}\label{sec3:eq2}
E(H_{d,j}^{s,l})\coloneqq\{\{\b,\g\}\mid\b+\g+\supp(p_j)\cap\B_d^{s,l}\ne\emptyset\}
\end{equation}
for $j=0,1,\ldots,m$, and further let
\begin{equation}\label{sec3:eq3}
\B_d^{s,l+1}=\bigcup_{j=0}^m\bigl(\supp(p_j)+\supp((H_{d,j}^{s,l})')\bigr).
\end{equation}
Doing so, we also get a finite ascending chain of support sets
$\B_d^{s,1}\subseteq\cdots\subseteq\B_d^{s,\tilde{l}}=\B_d^{s,\tilde{l}+1}=\cdots$
and a finite ascending chain of graphs
$(H_{d,j}^{s,1})'\subseteq\cdots\subseteq(H_{d,j}^{s,\tilde{l}})'=(H_{d,j}^{s,\tilde{l}+1})'=\cdots$
for each $j=0,1,\ldots,m$.

\begin{remark}
To handle polynomial inequalities imposed over another constraint set (say, defined by $q_j\ge0,j=1,\ldots,t$), we simply replace $p_j$ by $q_j$ in \eqref{sec3:eq2} and \eqref{sec3:eq3} to obtain the related sparsity graphs.
\end{remark}

\begin{example}\label{sec3:ex2}
Continue considering Example \ref{sec3:ex1}. Take $s=1$. The sparsity graph $H_{2,0}^{1,1}$ is displayed in Fig.~\ref{supp2}. If we use maximal chordal extensions, then $(H_{2,0}^{1,1})'$ consists of two cliques respectively of size $6,4$, and $\B_2^{1,2}=\A_{2}^2$. In fact, we have $\B_2^{1,2}=\B_2^{1,3}=\cdots$ and $(H_{2,j}^{1,2})'=(H_{2,j}^{1,3})'=\cdots$ for $j=0,1,2,3$.

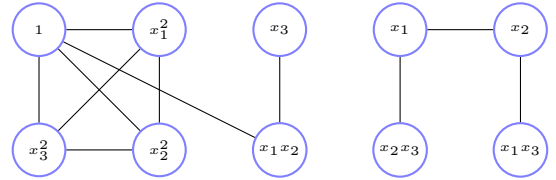
\begin{figure}[htbp]
\centering
{\tiny
\begin{tikzpicture}[scale=0.8, every node/.style={circle, draw=blue!50, thick, minimum size=7mm}]
\node (n1) at (-2,0) {$x_3$};
\node (n2) at (0,0) {$x_1$};
\node (n3) at (2,0) {$x_2$};
\node (n4) at (-2,-2) {$x_1x_2$};
\node (n5) at (0,-2) {$x_2x_3$};
\node (n6) at (2,-2) {$x_1x_3$};
\node (n7) at (-6,0) {$1$};
\node (n8) at (-4,0) {$x_1^2$};
\node (n9) at (-6,-2) {$x_3^2$};
\node (n10) at (-4,-2) {$x_2^2$};
\draw (n1)--(n4);
\draw (n7)--(n8);
\draw (n7)--(n9);
\draw (n7)--(n10);
\draw (n8)--(n9);
\draw (n8)--(n10);
\draw (n9)--(n10);
\draw (n2)--(n5);
\draw (n3)--(n6);
\draw (n2)--(n3);
\draw (n4)--(n7);
\end{tikzpicture}}
\caption{\small The sparsity graph $H_{2,0}^{1,1}$ for Example \ref{sec3:ex2}.}\label{supp2}
\end{figure}

\end{example}

\begin{remark}
The reason for using two different support sets for $v$ and $w$ is that $v$ and $w$ are of different degrees when one fixes a relaxation order $d$ for \eqref{sos}.
\end{remark}

The two indices $s$ and $l$ are used to control the size of support sets of $v$ and $w$ respectively. For a pair $s,l\ge1$, we may thereby consider the following sparse SOS relaxation for approximating the MPI set:
\begin{equation}\label{ssos}
\theta_d^{s,l}\coloneqq\begin{cases}
    \inf\limits_{a_j,b_j,c_j,v,w}\quad&\int_{X}w(\x)\,\mathrm{d}\x\\
    \quad\,\,\,\text{s.t.}&v\in\R[\A_d^s\cap\N^n_{2d+1-d_f}],w\in\R[\B_d^{s,l}],\\
    &\beta v-\nabla v\cdot \mathbf{f}=a_0+\sum_{j=1}^ma_jp_j,\\
    &w=b_0+\sum_{j=1}^mb_jp_j,\\
    &w-v-1=c_0+\sum_{j=1}^mc_jp_j,\\
    &a_j\in\Sigma[(G_{d,j}^{s})'],b_j,c_j\in\Sigma[(H_{d,j}^{s,l})'],\\
    &\quad\quad\quad\quad\quad\quad\quad\quad j=0,1,\ldots,m.
\end{cases}
\end{equation}
Notice that the Gram matrix of any SOS involved in \eqref{ssos} admits a block decomposition because of Theorem \ref{sec2-thm}. Hence the corresponding SDP could be easier to solve.

\begin{proposition}\label{sec3:prop}
With the above notations, we have $\theta_{d+1}^{s,l}\le\theta_d^{s,l}$, $\theta_{d}^{s+1,l}\le\theta_d^{s,l}$, $\theta_{d}^{s,l+1}\le\theta_d^{s,l}$, and $\theta_d^{s,l}\ge\theta_d$ for $d\ge\max\,\{\lceil d_f/2\rceil,\lceil d_p/2\rceil\},s\ge1,l\ge1$.
\end{proposition}
\begin{IEEEproof}
First note that $\supp(G_1)\subseteq\supp(G_2)$, $\Sigma[G_1]\subseteq\Sigma[G_2]$ if $G_1$ and $G_2$ are two graphs with $G_1\subseteq G_2$.
As a result, the feasible set of the corresponding sparse SOS relaxation becomes larger when we increase $d$, $s$, or $l$, from which the first three inequalities then follow. The last inequality $\theta_d^{s,l}\ge\theta_d$ holds because the feasible set of the sparse SOS relaxation is a subset of the corresponding dense SOS relaxation.
\end{IEEEproof}

\vspace{1em}
\noindent{\bf Sign symmetry.}
The {\em sign symmetries} of the system \eqref{dyn} with the constraint set $X$ \eqref{cons} consist of all vectors $\br=(r_i)\in\Z_2^n\coloneqq\{0,1\}^n$ that satisfy
$f_i((-1)^{\br}\circ\x)=(-1)^{r_i} f_i(\x)$ for $i=1,\ldots,n$,  
and 
$p_j((-1)^{\br}\circ\x)=p_j(\x)$ for $j=1,\ldots,m$,
where $(-1)^{\br}\coloneqq((-1)^{r_1},\ldots,(-1)^{r_n})$. 
Given a set of sign symmetries $R\subseteq\Z_2^n$, we define $R^{\perp}\coloneqq\{\a\in\N^n\mid\br\cdot\a\equiv0\,(\textrm{mod }2),\,\forall\,\br\in R\}$. Then a polynomial $g$ is invariant under the sign symmetries $R$\footnote{That is, $g((-1)^{\br}\circ\x)=g(\x)$ for any $\br\in R$.} if and only if $\supp(g)\subseteq R^{\perp}$. The following theorem\footnote{A similar symmetry reduction already appeared in \cite{fantuzzi2020bounding} in the study of bounding extreme events in dynamical systems.} tells us that the SOS relaxation \eqref{sos} inherits the sign symmetries of the dynamical system \eqref{dyn}.
\begin{theorem}\label{ss}
Let $R$ be the set of sign symmetries of the system \eqref{dyn} with the constraint set $X$ \eqref{cons}. If we additionally impose the constraints that $\supp(v),\supp(w)\subseteq R^{\perp}$ and $\supp(a_j),\supp(b_j),\supp(c_j)\subseteq R^{\perp}$ for $j=0,1,\ldots,m$ in \eqref{sos}, the resulting program has the same optimum with \eqref{sos}.
\end{theorem}
\begin{IEEEproof}
Suppose that $v,w,\{a_j\}_{j=0}^m,\{b_j\}_{j=0}^m,\{c_j\}_{j=0}^m$ are an optimal solution to \eqref{sos}. We remove the terms of $v,w$ with exponents not belonging to $R^{\perp}$ from the expression of $v,w$ and denote the resulting polynomials by $\tilde{v},\tilde{w}$ respectively. Let $Q_j$ be a PSD Gram matrix of $a_j$ for any $j$ such that $a_j=(\x^{\N^n_{d-d_j}})^\intercal Q_j\x^{\N^n_{d-d_j}}$. We then define $\widetilde{Q}_j\in\S^{|\N^n_{d-d_j}|}$ by letting $[\widetilde{Q}_j]_{\b\g}\coloneqq[Q_j]_{\b\g}$ for $\b+\g\in R^{\perp}$ and putting zeros elsewhere,
and let $\tilde{a}_j\coloneqq(\x^{\N^n_{d-d_j}})^\intercal\widetilde{Q}_j\x^{\N^n_{d-d_j}}$. One can easily check that $\widetilde{Q}_j$ is block diagonal (after an appropriate permutation on rows and columns). So $Q_j\succeq0$ implies $\widetilde{Q}_j\succeq0$ and it follows that $\tilde{a}_j$ is an SOS polynomial. In a similar way, we define $\tilde{b}_j$,$\tilde{c}_j$ for $j=0,1,\ldots,m$ which are all SOS polynomials by a similar argument as for $\tilde{a}_j$. As we remove exactly the terms with exponents not belonging to $R^{\perp}$ from both sides of the equations in \eqref{sos}, $\tilde{v},\tilde{w},\{\tilde{a}_j\}_{j=0}^m,\{\tilde{b}_j\}_{j=0}^m,\{\tilde{c}_j\}_{j=0}^m$ are again a feasible solution to \eqref{sos}. It remains to show $\int_{X}\tilde{w}(\x)\,\mathrm{d}\lambda=\int_{X}w(\x)\,\mathrm{d}\lambda$. Take any $\a\in\N^n_{2d}\setminus R^{\perp}$. Then there exists $\br=(r_i)_i\in R$ such that $\a\cdot\br\not\equiv0\,(\textrm{mod }2)$. We have
$\int_{X}\x^{\a}\,\mathrm{d}\lambda=\int_{X}((-1)^{\br}\circ\x)^{\a}\,\mathrm{d}\lambda=-\int_{X}\x^{\a}\,\mathrm{d}\lambda$,
where the first equality follows from the fact that $X$ is invariant under the sign symmetry $\br$. This immediately gives $\int_{X}\x^{\a}\,\mathrm{d}\lambda=0$ from which it follows $\int_{X}\tilde{w}(\x)\,\mathrm{d}\lambda=\int_{X}w(\x)\,\mathrm{d}\lambda$.
\end{IEEEproof}
From the proof of Theorem \ref{ss}, we see that the sign symmetries of the dynamical system \eqref{dyn} endow any SOS polynomial involved in \eqref{sos} with a block structure.
Our iterative procedure to exploit term sparsity actually produces block structures that are compatible with the sign symmetries of the dynamical system. Furthermore, when maximal chordal extensions are used in the construction, the block structures converge to the one given by the sign symmetries of the system. The key observation is the following lemma. 
\begin{lemma}\label{ss-lm}
Given $d\ge\max\,\{\lceil d_f/2\rceil,\lceil d_p/2\rceil\}$, the sign symmetries of the system \eqref{dyn} with the constraint set $X$ \eqref{cons} coincide with the sign symmetries of $\A_d^1$, i.e., the set $\{\br\in\Z_2^n\mid\br\cdot\a\equiv0\textrm{ mod }2,\,\forall\,\a\in\A_d^1\}$.
\end{lemma}
\begin{IEEEproof}
Let us denote the set of sign symmetries of the system \eqref{dyn} by $R$ and the set of sign symmetries of $\A_d^1$ by $R'$. For any $\br\in R$, we wish to show $\br\in R'$. It suffices to show $\br\cdot\a\equiv0\textrm{ mod }2$ for any $\a\in\A_d^1=\A\cup\supp(\nabla v\cdot \mathbf{f})\cup2\N^n_d$ where $\A=\bigcup_{j=1}^m\supp(p_j)$ and $v$ is a polynomial with generic coefficients supported on $\A$. If $\a\in\A\cup2\N^n_d$, we clearly have $\br\cdot\a\equiv0\textrm{ mod }2$ by definition. Noting $\nabla v(\x)=\nabla v((-1)^{\br}\circ\x)=(-1)^{\br}\circ(\nabla v)((-1)^{\br}\circ\x)$, we have $(\nabla v\cdot \mathbf{f})((-1)^{\br}\circ\x)=(\nabla v)((-1)^{\br}\circ\x)\cdot \mathbf{f}((-1)^{\br}\circ\x)=((-1)^{\br}\circ\nabla v(\x))\cdot ((-1)^{\br}\circ \mathbf{f}(\x))=(\nabla v\cdot \mathbf{f})(\x)$. 
As a result, $\br\cdot\a\equiv0\textrm{ mod }2$ for any $\a\in\supp(\nabla v\cdot \mathbf{f})$. 
Conversely, for any $\br\in R'$, we wish to show $\br\in R$. Since $\A\subseteq\A_d^1$, we have $p_j((-1)^{\br}\circ\x)=p_j(\x)$ for $j=1,\ldots,m$. Let $\a=(\alpha_i)\in\A$ with $\alpha_1>0$. Then $\supp(\x^{\a}f_1/x_1)\subseteq\A_d^1$. The condition $\br\in R'$ implies $(-1)^{\br\cdot\a}\x^{\a}f_1((-1)^{\br}\circ\x)/((-1)^{r_1}x_1)=\x^{\a}f_1/x_1$, which gives $f_1((-1)^{\br}\circ\x)=(-1)^{r_1}f_1(\x)$ as $(-1)^{\br\cdot\a}=1$. Similarly, we can prove $f_i((-1)^{\br}\circ\x)=(-1)^{r_i}f_i(\x)$ for $i=2,\ldots,n$. Thus $\br\in R$.
\end{IEEEproof}

Based on Lemma \ref{ss-lm}, we can prove the following theorem by a similar argument as for \cite[Theorem 6.5]{tssos} and so we omit the proof.
\begin{theorem}\label{thm:SignSym}
Given $d\ge\max\,\{\lceil d_f/2\rceil,\lceil d_p/2\rceil\}$, if maximal chordal extensions are used in the construction, then the block structures produced by the iterative procedure stated above Proposition \ref{sec3:prop} converge to the one given by the sign symmetries of the system \eqref{dyn} (with the constraint set $X$ \eqref{cons}) as $s,l$ increase. As a corollary of this fact and Theorem \ref{ss}, in this case, we have $\theta_d^{s,l}=\theta_d$ when $s,l$ are sufficiently large\footnote{The values of $s,l$ for $\theta_d^{s,l}=\theta_d$ to be valid are a priori unknown, but they are typically no greater than $3$ in practice.}.
\end{theorem}

\begin{remark}
Except maximal chordal extensions, other choices of chordal extensions may result in moment-SOS hierarchies with lower complexity, but that generally produce more conservative estimates for the problem at hand.
\end{remark}

\begin{remark}
From a more general point of view, the core elements of the strategy for exploiting term sparsity consist of (i) an algorithm to construct an increasing sequence $G_1 \subseteq G_2 \subseteq \cdots$ of sparsity graphs for the Gram matrices of SOS polynomials, and (ii) a way to impose the positive semidefiniteness of a large Gram matrix via decomposition into smaller matrix constraints. The strategy we suggest in this paper is just one particular choice (perhaps the most natural one as it converges to sign symmetries automatically).
\end{remark}

\begin{figure*}[!htb]
\centering
 \begin{minipage}{0.22\linewidth}
\includegraphics[scale=0.26]{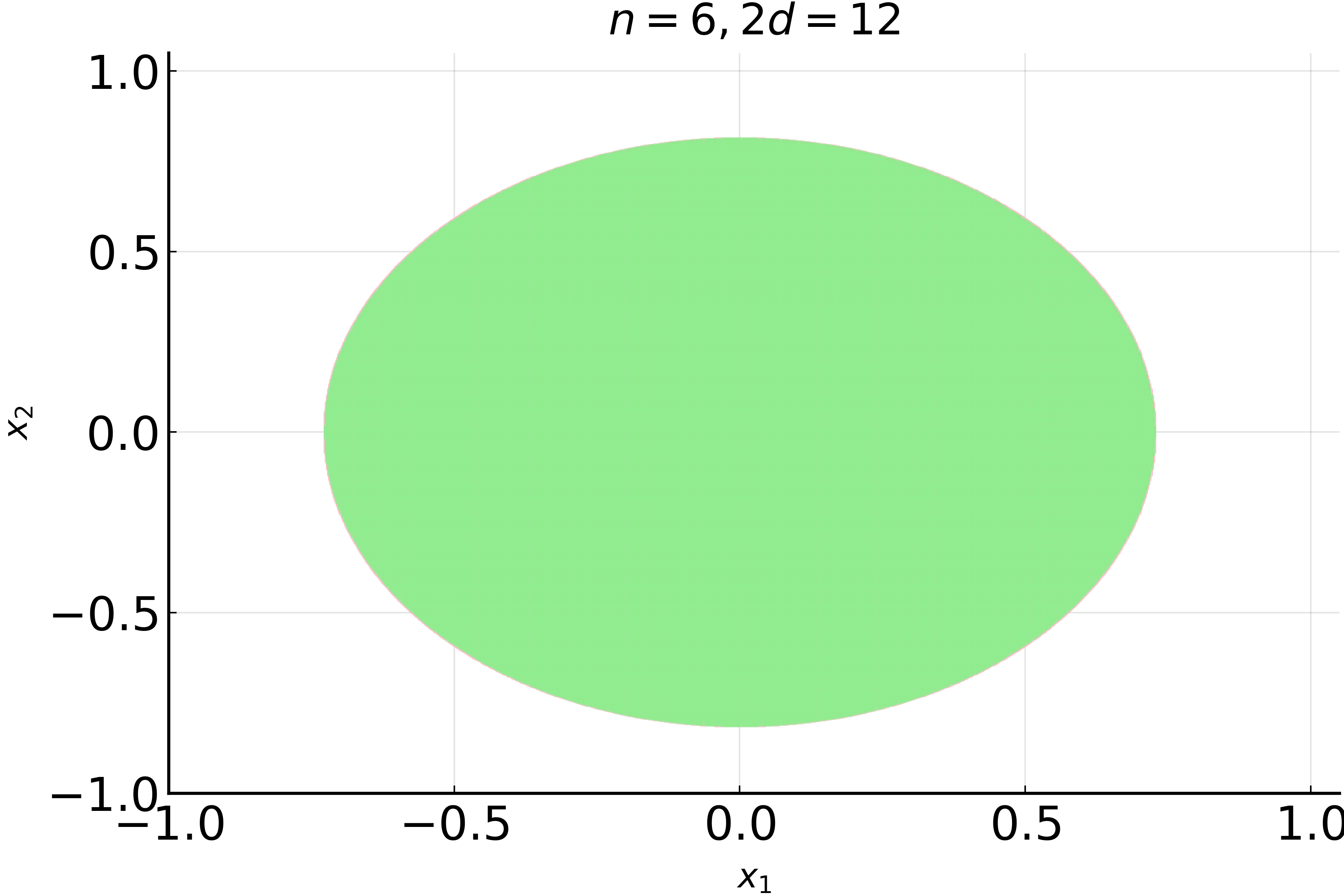}
\end{minipage}
\begin{minipage}{0.22\linewidth}
\includegraphics[scale=0.26]{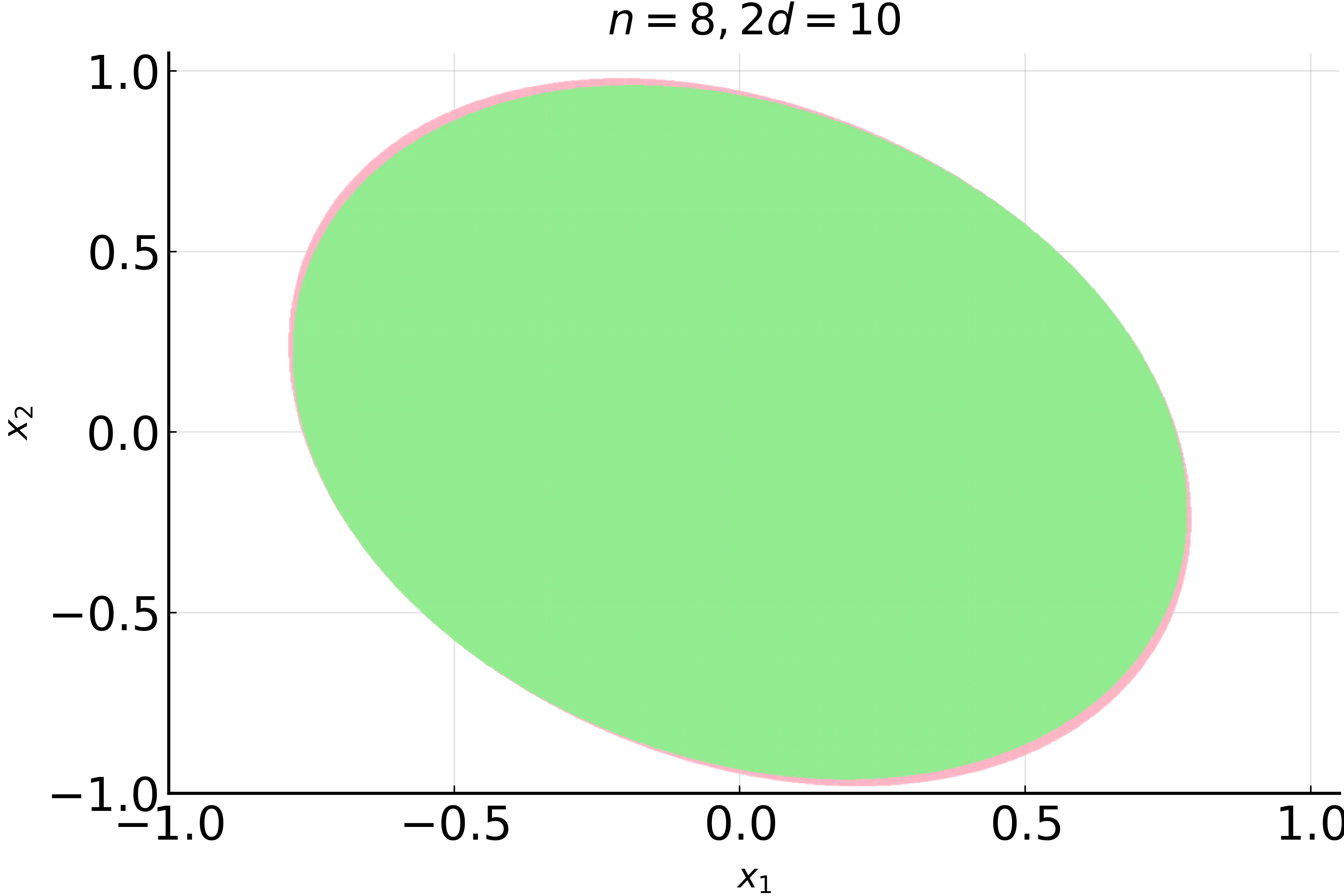}
\end{minipage}
\begin{minipage}{0.22\linewidth}
\includegraphics[scale=0.26]{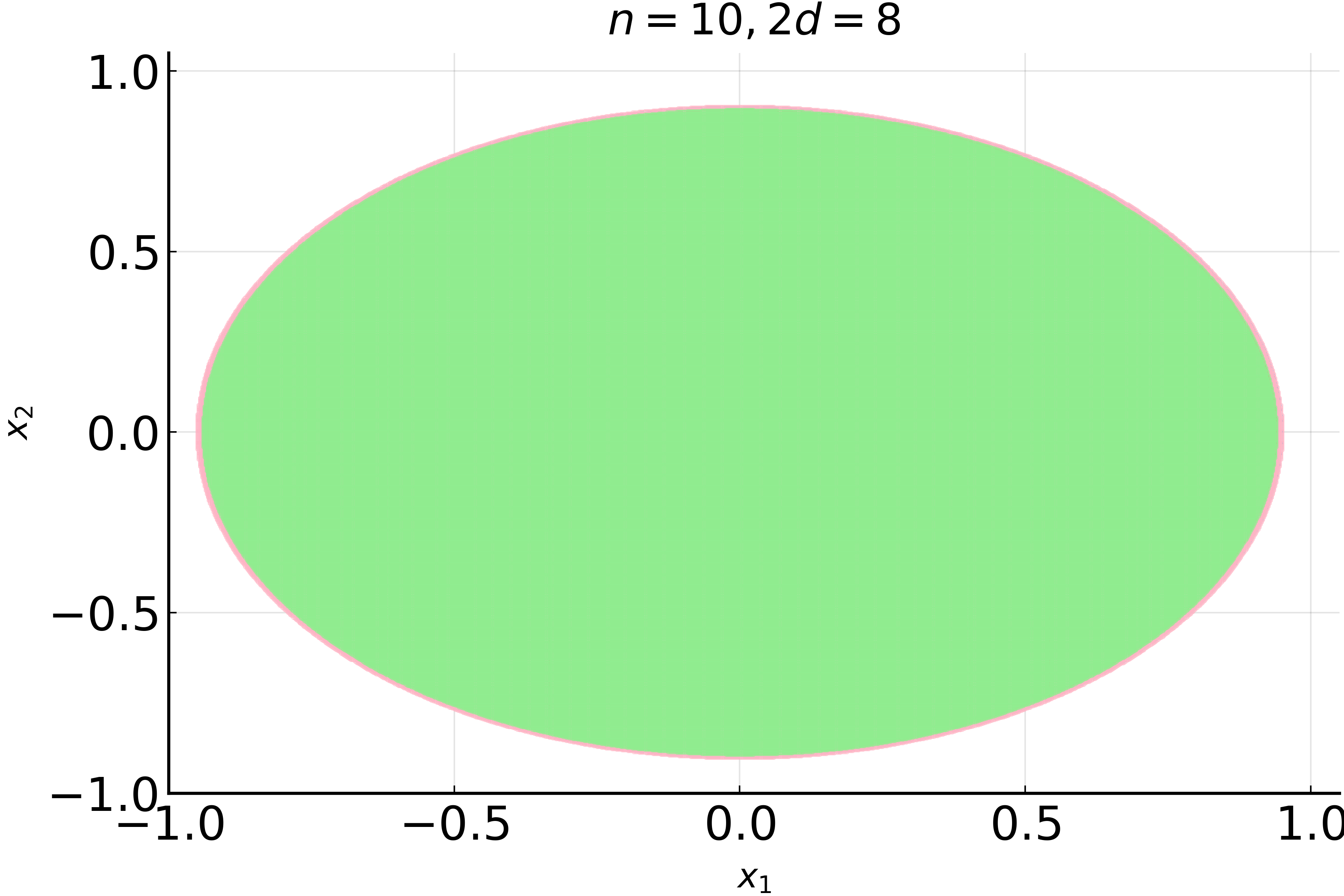}
\end{minipage}
\begin{minipage}{0.22\linewidth}
\includegraphics[scale=0.26]{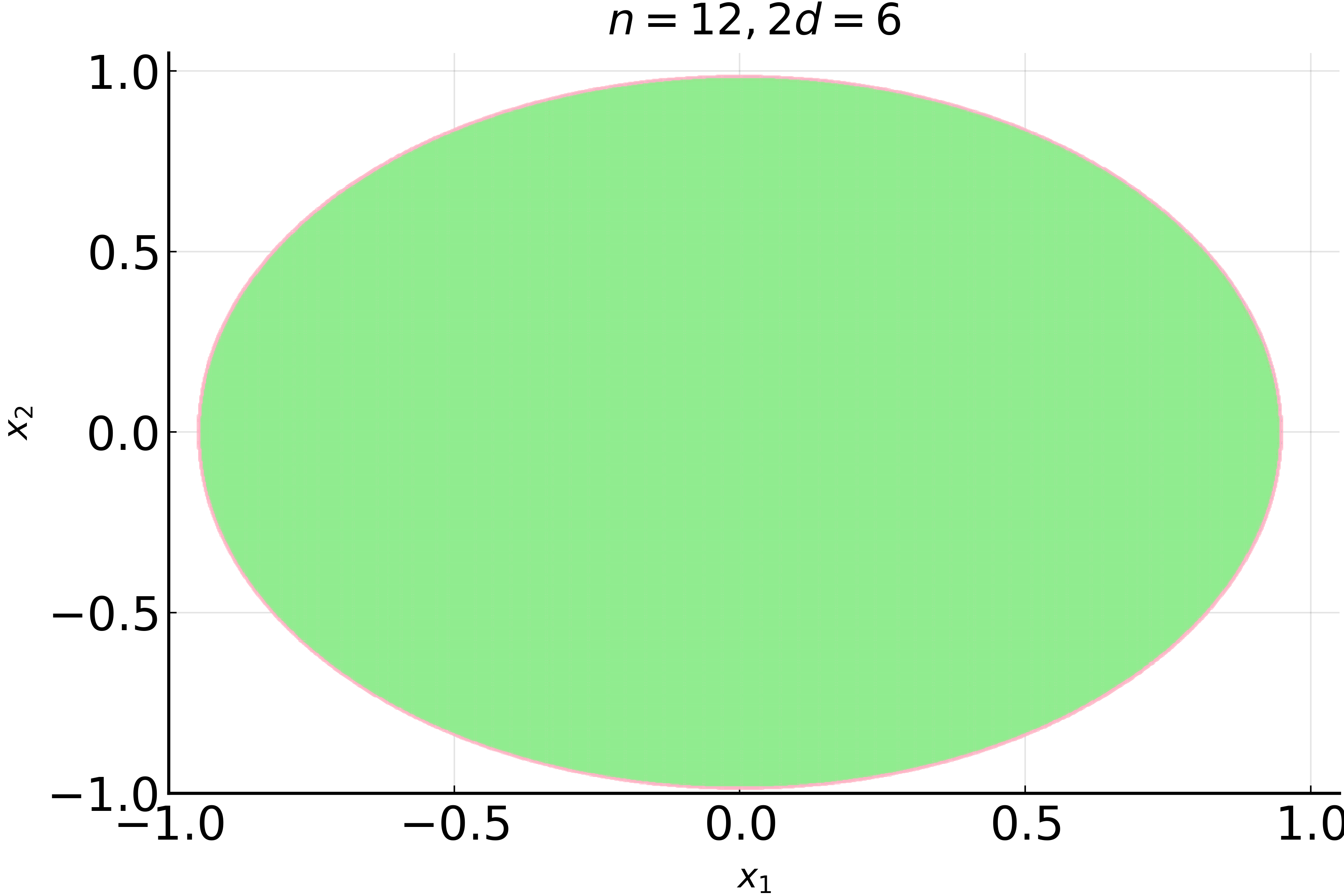}
\end{minipage}
\caption{\small Outer approximations for the MPI set of randomly generated models. TS {\scalebox{1.2}{\color{pink1} $\bullet$}}\; SS/FD \scalebox{1.2}{{\color{lightgreen} $\bullet$}}.}\label{random}
\end{figure*}

\section{Illustrative examples}\label{example}
In this section, we give several numerical examples to illustrate the iterative procedure to exploit term sparsity. The procedure has been implemented as a Julia package {\tt SparseDynamicSystem} which is available at:
\begin{center}
\url{https://github.com/wangjie212/SparseDynamicSystem}.
\end{center}
All examples were computed on an Intel i9-10900@2.80GHz CPU with 64GB RAM memory and {\tt Mosek} is used as an SDP solver. The timing is recorded in seconds. ``opt'' denotes the optimum. The index $l$ is fixed to $1$.

\vspace{-0.2cm}

\subsection{Comparison with the method of \cite{tacchi2020approximating}}\label{secA}
Consider $\{\dot{x_1}=(x_1^2+x_2^2-\frac{1}{4})x_1,\dot{x_2}=(x_2^2+x_3^2-\frac{1}{4})x_2,\dot{x_3}=(x_2^2+x_3^2-\frac{1}{4})x_3\}$
with the constraint set $X=[-1,1]^3$. According to \cite{tacchi2020approximating}, the system is decoupled w.r.t. $\{x_1,x_2\}$ and $\{x_2,x_3\}$. This system cannot be decoupled by \cite{schlosser2020sparse}. We compare the results for approximating the MPI set obtained by exploiting term sparsity ($s=1$), sign symmetries ($s=2$) \footnote{Namely, the block structures of the sparse hierarchy converge to the one determined by the sign symmetries of the system when $s=2$.}, the dense method and the method of \cite{tacchi2020approximating} in Table \ref{ex1}. We see that: TS is the fastest and provides medium bounds; SS and FD provide the best bounds whereas the method of \cite{tacchi2020approximating} provides the worst bounds.

\begin{table}[htbp]
\caption{\small The results for Section~\ref{secA}, TS: term sparsity, SS: sign symmetry, FD: fully dense.}\label{ex1}
\centering
\begin{tabular}{c|cc|cc|cc|cc}
\multirow{2}*{$2d$}&\multicolumn{2}{c|}{TS}&\multicolumn{2}{c|}{SS}&\multicolumn{2}{c|}{FD}&\multicolumn{2}{c}{\cite{tacchi2020approximating}}\\
\cline{2-9}
&opt&time&opt&time&opt&time&opt&time\\
\hline
18&2.86&0.15&1.66&2.55&1.66&33.0&6.21&0.57\\
20&2.86&0.21&1.55&3.72&1.55&56.2&5.89&1.31\\
22&2.86&0.33&1.49&6.00&1.49&132&5.64&2.38\\
\end{tabular}
\end{table}

\subsection{Comparison with the method of \cite{schlosser2020sparse}}\label{secB}
Consider the system $\{\dot{x_1}=(x_1^2+x_2^2-\frac{1}{4})x_1,
\dot{x_2}=x_2,
\dot{x_3}=(x_2^2+x_3^2-\frac{1}{4})x_3\}$
with the constraint set $X=[-1,1]^3$. According to \cite{schlosser2020sparse}, the system is decoupled w.r.t. $\{x_1,x_2\}$ and $\{x_2,x_3\}$. This system cannot be decoupled by \cite{tacchi2020approximating}. We compare the results for approximating the MPI set obtained by exploiting term sparsity ($s=1$), sign symmetries ($s=2$), the dense method, and the method of \cite{schlosser2020sparse} in Table \ref{ex2}. We can conclude that: SS and FD provide  the same bounds while the former is more efficient; TS and the method of \cite{tacchi2020approximating} are even more efficient but provide weaker bounds.

\begin{table}[htbp]
\caption{\small The results for Section~\ref{secB}, TS: term sparsity, SS: sign symmetry, FD: fully dense.}\label{ex2}
\centering
\begin{tabular}{c|cc|cc|cc|cc}
\multirow{2}*{$2d$}&\multicolumn{2}{c|}{TS}&\multicolumn{2}{c|}{SS}&\multicolumn{2}{c|}{FD}&\multicolumn{2}{c}{\cite{schlosser2020sparse}}\\
\cline{2-9}
&opt&time&opt&time&opt&time&opt&time\\
\hline
18&1.79&0.16&0.68&2.81&0.68&36.2&2.11&1.81\\
20&1.79&0.24&0.56&4.97&0.56&46.8&1.88&1.37\\
22&1.79&0.36&0.51&5.43&0.51&141&1.73&2.69\\
\end{tabular}
\end{table}

\vspace{-0.2cm}

\subsection{The extended Lorenz system}\label{secC}
Consider the extended Lorenz system $\{
    \dot{x_1}=10x_1-12x_2,
    \dot{x_2}=-\frac{70}{3}x_1+x_2+\frac{125}{3}x_1x_3,
    \dot{x_3}=\frac{8}{3}x_3-15x_1x_2,
    \dot{x_4}=10(x_4-x_1),
    \dot{x_5}=x_1(28-x_3)-x_5\}$
with the constraint set $X=[-1,1]^5$. According to \cite{schlosser2020sparse}, the system is decoupled w.r.t. $\{x_1,x_2,x_3,x_4\}$ and $\{x_1,x_2,x_3,x_5\}$. This system cannot be decoupled by \cite{tacchi2020approximating}. In Table \ref{ex3}, we report the results for approximating the MPI set by exploiting term sparsity ($s=2$), sign symmetries ($s=3$), the dense method and the method of \cite{schlosser2020sparse}. Here we see that: SS is several times faster than FD without sacrificing any accuracy; TS is slightly faster than SS while providing somewhat weaker bounds; the method of \cite{schlosser2020sparse} is the fastest but provides weaker bounds.

\begin{table}[htbp]
\caption{\small The results for Section~\ref{secC}, TS: term sparsity, SS: sign symmetry, FD: fully dense.}\label{ex3}
\centering
\begin{tabular}{c|cc|cc|cc|cc}
\multirow{2}*{$2d$}&\multicolumn{2}{c|}{TS}&\multicolumn{2}{c|}{SS}&\multicolumn{2}{c}{FD}&\multicolumn{2}{c}{\cite{schlosser2020sparse}}\\
\cline{2-9}
&opt&time&opt&time&opt&time&opt&time\\
\hline
8&3.31&1.37&3.24&2.29&3.24&3.48&4.20&0.62\\
10&2.59&7.64&2.45&8.72&2.45&35.3&3.19&2.76\\
12&1.59&50.3&1.53&60.1&1.53&271&1.86&14.0\\
\end{tabular}
\end{table}

\vspace{-0.2cm}

\subsection{Randomly generated models (\cite{tan2008stability})}
We consider the following sparse dynamical system for varying $n$:
$\dot{x}_i=(\x^\intercal B\x-1)x_i,i=1,\ldots,n$,
where $B\in\S_{G}$ is a random positive definite matrix satisfying:
\begin{enumerate}
    \item $G$ is a random graph with $n$ nodes and $n-4$ edges;
    \item For $1\le i\le n$, $B_{ii}\in[1,2]$ and for $1\le i<j\le n$, $B_{ij}\in[-0.5,0.5]$.
\end{enumerate}
The constraint set is $X=[-1,1]^n$. For each system, we approximate the MPI set by exploiting term sparsity ($s=1$), sign symmetries ($s=2$), and the dense method.
Fig.~\ref{random} shows outer approximations of the MPI set for $n=6,8,10,12$, respectively.
In Table \ref{random-t}, we list optima and running time for solving corresponding SDPs. From Fig.~\ref{random} and Table~\ref{random-t}, we can conclude: SS is faster than FD by one or two orders of magnitude without sacrificing any accuracy; TS is several times faster than SS while providing slightly weaker bounds.

\begin{table}[htb]
\caption{\small The results for randomly generated models, TS: term sparsity, SS: sign symmetry, FD: fully dense.}\label{random-t}
\centering
\begin{tabular}{c|c|cc|cc|cc}
\multirow{2}*{$n$}&\multirow{2}*{$2d$}&\multicolumn{2}{c|}{TS}&\multicolumn{2}{c|}{SS}&\multicolumn{2}{c}{FD}\\
\cline{3-8}
&&opt&time&opt&time&opt&time\\
\hline
\multirow{2}*{$6$}&10&6.60&5.05&6.23&19.1&6.23&1516\\
&12&5.56&118&5.05&183&-&-\\
\hline
\multirow{2}*{$8$}&8&13.9&28.8&13.4&46.6&13.4&1909\\
&10&11.4&281&10.8&1196&-&-\\
\hline
\multirow{2}*{$10$}&6&93.3&5.58&93.3&11.8&93.3&359\\
&8&30.7&528&29.8&586&-&-\\
\hline
\multirow{2}*{$12$}&4&510&0.37&438&1.02&438&5.63\\
&6&312&76.5&311&130&311&2272\\
\end{tabular}
\end{table}

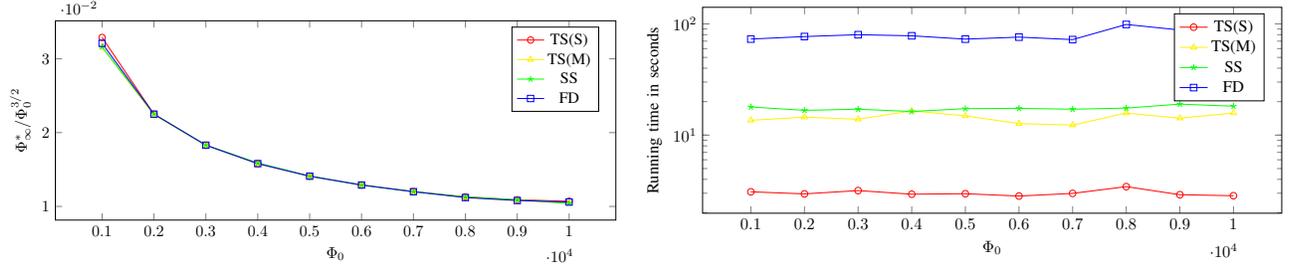
\begin{figure*}[!htb]
\centering
\begin{minipage}{0.46\linewidth}
\begin{tikzpicture}[scale=0.6]
\begin{axis}[xlabel={$\Phi_0$},
ylabel={$\Phi^{\ast}_{\infty}/\Phi_0^{3/2}$},
legend pos=north east]
\addplot[color=red,mark=o]
coordinates{(1000,0.0329)(2000,0.0225)(3000,0.0183)(4000,0.0158)(5000,0.0141)(6000,0.0129)(7000,0.0120)(8000,0.0113)(9000,0.0109)(10000,0.0107)};
\addlegendentry{TS(S)}
\addplot[color=yellow,mark=triangle]
coordinates{(1000,0.0316)(2000,0.0224)(3000,0.0183)(4000,0.0158)(5000,0.0141)(6000,0.0129)(7000,0.0120)(8000,0.0112)(9000,0.0108)(10000,0.0105)};
\addlegendentry{TS(M)}
\addplot[color=green,mark=star]
coordinates{(1000,0.0316)(2000,0.0225)(3000,0.0183)(4000,0.0159)(5000,0.0141)(6000,0.0129)(7000,0.0120)(8000,0.0113)(9000,0.0109)(10000,0.0104)};
\addlegendentry{SS}
\addplot[color=blue,mark=square]
coordinates{(1000,0.0321)(2000,0.0225)(3000,0.0183)(4000,0.0158)(5000,0.0141)(6000,0.0129)(7000,0.0120)(8000,0.0112)(9000,0.0108)(10000,0.0106)};
\addlegendentry{FD}
\end{axis}
\end{tikzpicture}
\end{minipage}
\begin{minipage}{0.46\linewidth}
\begin{tikzpicture}[scale=0.62]
\begin{axis}[xlabel={$\Phi_0$}, ylabel={Running time in seconds}, legend pos=north east, ymode=log, log basis y=10]
\addplot[color=red,mark=o]
coordinates{(1000,3.09)(2000,2.96)(3000,3.17)(4000,2.94)(5000,2.97)(6000,2.83)(7000,2.99)(8000,3.44)(9000,2.91)(10000,2.85)};
\addlegendentry{TS(S)}
\addplot[color=yellow,mark=triangle]
coordinates{(1000,13.6)(2000,14.5)(3000,13.9)(4000,16.5)(5000,14.9)(6000,12.7)(7000,12.3)(8000,15.8)(9000,14.2)(10000,15.8)};
\addlegendentry{TS(M)}
\addplot[color=green,mark=star]
coordinates{(1000,17.9)(2000,16.7)(3000,17.1)(4000,16.3)(5000,17.3)(6000,17.4)(7000,17.1)(8000,17.5)(9000,19.0)(10000,18.2)};
\addlegendentry{SS}
\addplot[color=blue,mark=square]
coordinates{(1000,73.1)(2000,77.0)(3000,80.0)(4000,78.2)(5000,73.0)(6000,76.1)(7000,72.3)(8000,98.9)(9000,88.2)(10000,74.2)};
\addlegendentry{FD}
\end{axis}
\end{tikzpicture}
\end{minipage}
\caption{\small Upper bounds on $\Phi^{\ast}_{\infty}$ and running time for the 16 mode fluid model with $2d=4$. Here, the values of $\Phi^{\ast}_{\infty}$ are normalized by $\Phi_0^{3/2}$ to be better visualized.}\label{fg1}
\end{figure*}

\subsection{Bounding extreme events for a 16 mode fluid model}
This model (\cite[Example 4.2]{fantuzzi2020bounding}) is given by the system:
\begin{equation}\label{16model}
\dot{x}_n=-(2\pi n)^2x_n+\sqrt{2}\pi n\left(\sum_{i=1}^{16-n}x_ix_{i+n}-\frac{1}{2}\sum_{i=1}^{n-1}x_ix_{n-i}\right),
\end{equation}
for $n=1,\ldots,16$.
Let $\Phi(\x)=2\pi^2\sum_{i=1}^{16}i^2x_i^2$, $X_0=\{\x\in\R^{16}\mid\Phi(\x)=\Phi_0\}$ and $X=\{\x\in\R^{16}\mid||\x||_2^2\le\Phi_0/(2\pi^2)\}$, where $\Phi_0\in\R$.
Let $\Phi^{\ast}_{\infty}$ denote the largest value attained by $\Phi(\x(t;t_0,\x_0))$ among all trajectories that start from $X_0\subseteq X$ and evolve forward over the time
interval $[0,\infty)$. SOS relaxations were proposed in \cite{fantuzzi2020bounding} to bound $\Phi^{\ast}_{\infty}$ from above. This system cannot be decoupled by either \cite{schlosser2020sparse} or \cite{tacchi2020approximating}. We can adapt the iterative procedure presented in this paper to derive sparse SOS relaxations for bounding $\Phi^{\ast}_{\infty}$. Taking the relaxation order $d=2$, we solve the dense relaxation (FD) and the sparse relaxations with $s=1$ and $s=2$ (SS). In particular, for $s=1$, we include the results obtained either with approximately smallest chordal extensions (TS(S)) or with maximal chordal extensions (TS(M)). The results are shown in Fig.~\ref{fg1}. 
It turns out that the upper bounds given by the four methods agree to within $4$\% for different $\Phi_0$. On average, TS(S) is six times faster than SS, TS(M) is slightly faster than SS, and SS is five times faster than FD, which indicates that there is significant speed-up by exploiting term sparsity.



\section{Conclusions}\label{con}
This paper presents a reduction approach by exploiting term sparsity for
the moment-SOS hierarchy of problems arising from the study of dynamical systems. As demonstrated by the numerical examples, this approach provides a trade-off between computational costs and the solution accuracy. Moreover, it is able to guarantee convergence under certain conditions and recover the sign symmetry reduction.

\bibliographystyle{IEEEtran}
\bibliography{refer}

\end{document}